# Symmetry Properties of Generalized Regular Polytopes

*Alexander V. Kharchenko*

*Email: alex@newhavensoft.net*

A concept of generalized regular polytope is introduced in this work. The number of its $(1...n-1)$-dimensional elements is not necessarily integer, though all the combinatorial and metric properties meet those of regular polytopes in a classic sense. New relationships between Schlafli symbol $\{f_1,...,f_{n-1}\}$ of the regular polytope and its metric parameters have been established. Using the generalized regular polytopes concept, group and metric properties of arbitrary metric space tessellations into regular honeycombs were investigated. It has been shown that sequential tessellations of space into regular honeycombs determine an infinite discrete group, having finite cyclic, dihedral, symmetric, and other subgroups. Set of generators and generating relations of the group are identified. Eigenvectors of regular honeycombs have been studied, and some of them shown to correspond to Schlafli symbols of known integer regular polytopes in 3 and 4 dimensions. It was discovered that group of all regular honeycombs comprises subsets having eigenvectors inducing a metric of the $(p, q)$ signature, and in particular, $(+---)$. These eigenvectors can be interpreted as self-reproducing generalized regular polytopes (eigentopes).

# Legend

| For 3-dimensional polytopes || For 4-dimensional polytopes ||
| Symbol | Description | Symbol | Description |
|---|---|---|---|
| **y** | number of vertices | **Y** | number of vertices |
| **x** | number of edges | **X** | number of edges |
| **n** | number of faces | **N** | number of faces |
|  |  | **Z** | number of 3-dimensional polytopes |
| **i** | index of vertex | **I** | index of vertex |
| **m** | number of edges incident to a face | **M** | number of edges incident to a vertex |
|  |  | **H** | number of 3-dimensional polytopes incident to a vertex |
|  |  | **U** | number of 3-dimensional polytopes incident to an edge |
| **{m, i}** | symbol of polytope | **{m, i, U}** | symbol of polytope |

**Concept of generalized regular polytopes**

Let us try to figure a regular polygon with a number of edges and vertices **m** = 4.5 on a plane (Fig. 1). For this purpose, we fix coordinates of a polygon center O and one of its vertices A, so that radius of the circum-circle is equal to $R$.

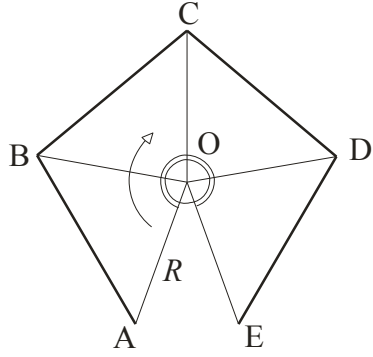

Fig. 1. An attempt to draw regular polygon with the number of edges and vertices **m** = 4.5

Then, moving clockwise from the vertex A, we will cut out consecutive angles $2\pi/\mathbf{m}$. Limiting rays by the radius $R$, we will obtain new vertices B, C, D and E. The last fifth vertex E appears to be not only redundant at 1/2, but also asymmetric with respect to the proximate vertices A and D, which contradicts with the initial requirement that polygon has to be *regular* one. We do not extend this construction to *stellated* polygons since their metric properties (e. g. perimeter, area) are not consistent with those of integer regular polygons.

Therefore non-integer regular polygon cannot be drawn in full neither on a Euclidean plane nor on any other surface despite of the fact that all metric relationships between its elements *meet* those of a conventional concept of regular polygon. It is possible however, fixing a center and one of its vertices, and taking into account equality of all the vertices, estimate uncertainty in an angle coordinate of all other vertices:

$$\Delta\varphi = \pi\frac{|\Delta\mathbf{m}|}{\mathbf{m}} ,  \qquad (1)$$

where $\Delta\mathbf{m}$ - the deviation of **m** from the nearest integer value.

Further we will consider an example of regular tessellation of a 3-dimensional space, which allows interpreting non-integer regular polytopes as *statistical* honeycombs.

Leaving an attempt to fit non-integer regular polygon into plane, we *represent* it by three coplanar vectors $^2\mathbf{R}_0$, $^2\mathbf{R}_1$ and $^1\mathbf{R}_0$ (Fig. 2), such that:

$$^2\mathbf{R}_0 = {}^2\mathbf{R}_1 + {}^1\mathbf{R}_0 \; ; \quad \left({}^2\mathbf{R}_1, {}^1\mathbf{R}_0\right) = 0, \qquad (2)$$

out from which any two vectors are linearly independent, and completely define regular polygon {**m**}, since:

$$\cos\frac{\pi}{\mathbf{m}} = \sqrt{\frac{\left({}^2\mathbf{R}_1\right)^2}{\left({}^2\mathbf{R}_0\right)^2}} \quad \text{and} \quad R = \sqrt{\left({}^2\mathbf{R}_0\right)^2}. \qquad (3)$$

Thus, triangle system of vectors $^2\mathbf{R}_0$, $^2\mathbf{R}_1$ and $^1\mathbf{R}_0$ forms characteristic 2-dimensional simplex of a regular polygon.

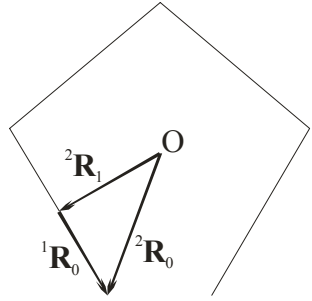

Fig. 2. Characteristic simplex of a regular polygon with the number of edges and vertices **m** = 4.5

For the *n*-dimensional regular polytope, characteristic simplex has been its *n*-dimensional simplex (orthoscheme). As Donald Coxeter showed [1], the characteristic simplex of regular polytope has simultaneously been a fundamental region of its symmetry group.

We have to abandon Berger's definition of a regular polytope based on the transitivity of its isotropy group on the set of its flags [2], as we concede that this set in general is not integer. Moreover, by the isotropy group of a polytope we are going to denote another set.

Another definition of an *abstract regular polytopes* is contained in the recent work of McMullen and Schulte [3], however it is similar to [2] in the part of transitivity of isotropy group on the integer set of flags, and therefore does not permit irrational numbers of all *n*-dimensional elements (faces).

Recurrent definitions of a regular polytope seem to be more appropriate; however they also don't bring sufficient clarity for the constructive use of metric relationships.

We will give an invariant definition of a generalized regular polytope in *n*-dimensional space with an arbitrary metric.

**Definition.** Generalized *n*-dimensional regular polytope is a polytope with characteristic simplex formed by system of $n(n+1)/2 - 1$ vectors $^a\mathbf{R}_b$ (orthoscheme), connecting centers of *a*-cells with centers of *b*-cells incident to them ($a > b$), such that:

$$^i\mathbf{R}_k = {}^i\mathbf{R}_j + {}^j\mathbf{R}_k, \quad 0 \le k < j < i \le n, \qquad (4)$$

and satisfying orthogonal relationships:

$$\left({}^{i}\mathbf{R}_{j}, {}^{k}\mathbf{R}_{l}\right) = 0, \quad 0 \leq l < k \leq j < i \leq n, \tag{5}$$

out from which any $n$ linearly independent vectors define a polytope unambiguously.

This definition obviously results in metric relationships and the structure of scalar product of constituting vectors:

$$\left({}^{i}\mathbf{R}_{k}\right)^{2} = \left({}^{i}\mathbf{R}_{j}\right)^{2} + \left({}^{j}\mathbf{R}_{k}\right)^{2}; \tag{6}$$

$$\left({}^{i}\mathbf{R}_{j}, {}^{i}\mathbf{R}_{k}\right) = \left({}^{i}\mathbf{R}_{j}\right)^{2}, \quad 0 \leq k < j < i \leq n. \tag{7}$$

Schlafli symbols $\{f_1, \ldots, f_{n-1}\}$ of the regular polytope may be directly expressed in terms of scalar squares of constituting vectors as follows:

$$\cos^{2}\frac{\pi}{f_{1}} = \frac{\left({}^{2}\mathbf{R}_{1}\right)^{2}}{\left({}^{2}\mathbf{R}_{0}\right)^{2}}; \tag{8}$$

$$\cos^{2}\frac{\pi}{f_{i}} = \frac{\left({}^{i+1}\mathbf{R}_{i}\right)^{2}}{\left({}^{i+1}\mathbf{R}_{i-1}\right)^{2}} \frac{\left({}^{i-1}\mathbf{R}_{i-2}\right)^{2}}{\left({}^{i}\mathbf{R}_{i-2}\right)^{2}}, \quad i = 2, \ldots, n-1. \tag{9}$$

On the other hand, these formulae allow to determine unambiguously metric relationships in a regular polytope with an arbitrary symbol, whose elements $\{f_1, \ldots, f_{n-1}\}$ are non-zero real or complex numbers.

It should be noted that proof of the formula equivalent to (8), is contained in [2]. The particular case of the formula (9) for the 3-dimensional regular polyhedra (i. e. for $f_2 := \mathbf{i}$) can be easily deduced from relationships obtained by Donald Coxeter in his classic work [1].

Hereinafter, terms «polytope», «regular polytope», and «regular honeycomb» will be used to denote «generalized regular polytope».

## Regular polytopes in 4-dimensional space

For the constituting vectors of 4-dimensional regular polytopes and their scalar squares the following special designations will be used:

$$\mathbf{p}_0 = {}^{4}\mathbf{R}_0 \quad \mathbf{p}_1 = {}^{4}\mathbf{R}_1 \quad \mathbf{p}_2 = {}^{4}\mathbf{R}_2 \quad \mathbf{p}_3 = {}^{4}\mathbf{R}_3$$

$$\rho_0 = \left(\mathbf{p}_0\right)^2 \quad \rho_1 = \left(\mathbf{p}_1\right)^2 \quad \rho_2 = \left(\mathbf{p}_2\right)^2 \quad \rho_3 = \left(\mathbf{p}_3\right)^2$$

**Relative basis in the space of polytopes**

Taking into account (6-9), symbol elements of the 4-dimensional regular polytope may be expressed by the following relationships:

$$\varepsilon \equiv \cos^2 \frac{\pi}{f_1} = \frac{\rho_1 - \rho_2}{\rho_0 - \rho_2}; \quad \delta \equiv \cos^2 \frac{\pi}{f_2} = \frac{(\rho_2 - \rho_3)(\rho_0 - \rho_1)}{(\rho_1 - \rho_3)(\rho_0 - \rho_2)}; \quad \eta \equiv \cos^2 \frac{\pi}{f_3} = \frac{\rho_3(\rho_1 - \rho_2)}{\rho_2(\rho_1 - \rho_3)}. \quad (10)$$

Variables $\varepsilon$, $\delta$, $\eta$ are defined as coordinates of the regular polytope in a relative basis $[E]$ of the space of regular polytopes whose vectors are triplets of various values $[\varepsilon, \delta, \eta]$. These will be denoted as $E$-symbols, or $E$-vectors of regular polytope due to their obvious concordance to the symbols $\{f_1, f_2, f_3\}$.

Also, we introduce another relative basis $[H]$ in the same space of polytopes:

$$\alpha = \frac{\rho_0}{\rho_3}, \quad \beta = \frac{\rho_0}{\rho_2}, \quad \gamma = \frac{\rho_0}{\rho_1}, \quad (11)$$

which is unambiguously associated with the basis $[E]$ by the relationships:

$$\alpha = \frac{(1-\varepsilon-\delta)(1-\delta-\eta)}{\varepsilon \, \delta \, \eta}, \quad \beta = \frac{(1-\varepsilon)(1-\delta-\eta)}{\varepsilon \, \delta}, \quad \gamma = \frac{1-\delta-\eta}{\varepsilon(1-\eta)}. \quad (12)$$

Respective vectors will be called $H$-symbols, or $H$-vectors of the regular polytope.

Hereinafter, relative bases $[E]$ and $[H]$ will be used interchangeably in various expressions depending on their simplicity. These bases are considered equivalent as long as corresponding Jacobians are not equal to zero.

**Basis of space, frame of polytope, and metric tensor**

Selection of basis vectors of space is limited by the requirement of their linear independence. From the entire set of possible bases, it is convenient to choose one, such that transformations of polytopes as a result of space tessellation into regular honeycombs have the simplest form. In author's opinion, set of vectors $\{\mathbf{p}_0, \mathbf{p}_1, \mathbf{p}_2, \mathbf{p}_3\}$ does have such properties. Further it will be called natural frame of polytope in the natural basis, and denoted as $\{\mathbf{P}\}$. Matrix of metric tensor in the natural basis $\{\mathbf{P}\}$ has the following form:

$$G_\mathbf{P} = \begin{pmatrix} \rho_0 & \rho_1 & \rho_2 & \rho_3 \\ \rho_1 & \rho_1 & \rho_2 & \rho_3 \\ \rho_2 & \rho_2 & \rho_2 & \rho_3 \\ \rho_3 & \rho_3 & \rho_3 & \rho_3 \end{pmatrix}. \quad (13)$$

Passing from natural to the orthogonal basis and the respective orthogonal frame $\{A\} = \{p_3, {}^3R_2, {}^2R_1, {}^1R_0\}$, we'll obtain diagonal matrix of metric tensor which is more convenient to analyze the metric's signature:

$$G_A = \rho_0 \begin{pmatrix} \frac{1}{\alpha} & & & 0 \\ & \frac{1}{\beta} - \frac{1}{\alpha} & & \\ & & \frac{1}{\gamma} - \frac{1}{\beta} & \\ 0 & & & 1 - \frac{1}{\gamma} \end{pmatrix}. \qquad (14)$$

The matrix of the coordinates change has been a symmetric matrix $B$, such that:

$$G_A = B\, G_P\, B, \quad \text{where} \quad B = \begin{pmatrix} 0 & 0 & 0 & 1 \\ 0 & 0 & 1 & -1 \\ 0 & 1 & -1 & 0 \\ 1 & -1 & 0 & 0 \end{pmatrix}. \qquad (15)$$

Note, that metric (14) is Euclidean one, provided that in the relative basis $[H]$ the following condition is satisfied:

$$\alpha > \beta > \gamma > 1. \qquad (16)$$

Otherwise, metric is pseudo-Euclidean one, and in particular if

$$\rho_0 > 0,\ 0 < \alpha < \beta < \gamma < 1, \qquad (17)$$

then metric (14) has a signature (+−−−).

In such a way symbol of the regular polytope (or corresponding vector in a relative basis) determines the metric of enclosing space with an accuracy of a constant factor. We do not use the orthonormal basis, as the majority of relationships are more complicated in it.

As an example, all the known integer regular polytopes may be presented:

| Polytope | $\alpha$ | $\beta$ | $\gamma$ |
|---|---|---|---|
| {3, 3, 3} | 16 | 6 | 8/3 |
| {4, 3, 3} | 4 | 2 | 4/3 |
| {3, 3, 4} | 4 | 3 | 2 |
| {3, 4, 3} | 2 | 3/2 | 4/3 |
| {5, 3, 3} | $4(7-3\sqrt{5})$ | $2(5-2\sqrt{5})$ | $4(1-\sqrt{5}/3)$ |
| {3, 3, 5} | $4(7-3\sqrt{5})$ | $(9/2)(1-3\sqrt{5})$ | $2(1-\sqrt{5}/5)$ |

All listed polytopes satisfy the condition (16), i.e. they are all correspond to the Euclidean metric. Marginal relationship $\alpha = \beta = \gamma = 1$ defines 3-dimensional sphere $S^3$ in a 4-dimensional Euclidean space.

Here are also examples of integer regular polytopes satisfying the condition (17), i.e. existing in a pseudo-Euclidean space of $(+---)$ signature:

| Polytope | $\alpha$ | $\beta$ | $\gamma$ |
|---|---|---|---|
| {5, 3, 4} | $7-3\sqrt{5}$ | $5-2\sqrt{5}$ | $3-\sqrt{5}$ |
| {4, 3, 5} | $7-3\sqrt{5}$ | $(3-\sqrt{5})/2$ | $1-\sqrt{5}/5$ |
| {3, 5, 3} | $4(9-4\sqrt{5})$ | $(3/2)(7-3\sqrt{5})$ | $2(1-\sqrt{5}/3)$ |
| {5, 3, 5} | $2(47-21\sqrt{5})$ | $(25-11\sqrt{5})/2$ | $4(1-2\sqrt{5}/5)$ |

Polytopes {4, 3, 5}, {3, 5, 3} and {5, 3, 5} were described first by Schlegel [13] as «regular hyperbolic honeycombs».

Hereinafter we will use an arbitrary metric of space keeping in mind a freedom in choosing symbol of a regular polytope.

### Tessellation of space into regular honeycombs

Only one regular tessellation of 3-dimensional space by cubes {4, 3}, and three regular tessellations of 4-dimensional space by polytopes {4, 3, 4}, {3, 4, 3}, and {3, 3, 4} are known. Condition of existence of such tessellations was formulated in [2] (in our designations) as follows:

«Tessellation of $n$-dimensional space into regular polytopes $\{f_1,..., f_{n-1}\}$ exists only if regular polytope $\{f_2,..., f_n\}$ does exist, and the regular polytope $\{f_1,..., f_n\}$ satisfy the relationship:

$$\frac{\left(^1\mathbf{R}_0\right)^2}{\left(^{n+1}\mathbf{R}_0\right)^2} = 0 \text{ ».} \tag{18}$$

In particular, for the cube {4, 3} there exists an infinite polytope {4, 3, 4} meeting the condition (18). This condition is equivalent to the relationship $\gamma = 1$ in the relative basis [$H$]. It should be noted that covering factor for the regular tessellation of 3-dimensional space by cubes is equal to 8, that is each node of tessellation belongs to 8 adjacent cubes. However, real tessellations in a physical 3-dimensional space tend to a rather minimal covering factor = 4, as for example in the uniform polycrystalline alloys and minerals having the minimal anisotropy between different crystallographic directions. Many similar «material» aspects of space division into regular statistical honeycombs, including one-dimensional case, were addressed in the work of Coxeter and Kharchenko [4].

Interesting results on the estimation of statistic average parameters of tessellation of 3-dimensional space with the minimal covering factor were obtained by Bernal [5] experimentally, by uniform compressing of chaotically placed and equally-sized plasticine balls. Count showed that resulting polyhedra have on the average 13.3 faces, and each face more frequently is bounded by 5 edges.

In works [6, 7], generalized thereafter by Miles [9], models of random division of space into polyhedral cells as a result of simultaneous growth of crystal nuclei with an equal rate in all directions, and uniformly distributed in $E_3$, were studied. They assumed that growth of any two crystals is terminated in a given point as they contact each other. For the numbers of vertices, edges, faces, and edges around each face the following average values were obtained: **y** = 27.07, **x** = 40.61, **n** = 15.54, and **m** = 5.23 respectively. Authors of [8] considered a similar model, but with a complementary condition that nuclei appear with the constant rate per unit volume. Taking this into account and using complicated statistical calculations, Meijering [6] obtained the following average values: **y** = 22.56, **x** = 33.84, **n** = 13.28, and **m** = 5.096. Another theoretical model suggested by Coxeter [10] yields: **y** = 23.13, **x** = 34.69, **n** = 13.56, and **m** = 5.115 .

The problem of «thirteen spheres» is also very close to the discussed matter. A new proof of the fact that maximum number of spheres contacting another given sphere of the same radius is equal to 13, has been recently presented by Anstreicher [11]. By removing central sphere and transforming other spheres to pyramids, we obtain a similar polyhedron with **n** = 13.

We will construct a regular tessellation of 3-dimensional space by reflecting the center of regular polytope {**m**, **i**} over its vertices (Fig. 3). Polyhedral angle of reflected polytope {**m**′, **i**′} is figured by the bold lines, and the direction to its center is marked by an arrow. Such a way of tessellation yields in a covering factor equal to the number of faces **n**, while the vertex index **i**′ of reflected polytope is equal to the number of edges **m** bounding each face of the initial polytope. To determine a value of **m**′, it is sufficient to express an angle $(\pi - \angle AOB)/2$ in terms of **m** and **i**, using relationships (6-9).

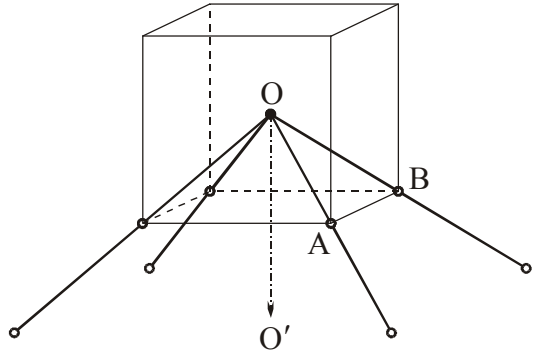

Fig. 3. Reflection of the center of 3-dimensional polytope over vertices

As a result of simple transformations, for the polytope $\{\mathbf{m}', \mathbf{i}'\}$ we will obtain:

$$\begin{cases} \cos^2 \dfrac{\pi}{\mathbf{m}'} = 1 - \dfrac{\cos^2 \dfrac{\pi}{\mathbf{m}}}{1 - \cos^2 \dfrac{\pi}{\mathbf{i}}} \; ; \\ \mathbf{i}' = \mathbf{m} \; . \end{cases} \quad (19)$$

These relationships have the following form in a relative basis $[E]$:

$$<\mathbf{A}>_{E3} = \begin{cases} \varepsilon' = 1 - \dfrac{\varepsilon}{1-\delta} \; ; \\ \delta' = \varepsilon \; , \end{cases} \quad (20)$$

where index $E3$ denotes basis of transformations and dimension of space.

With the transformations $<\mathbf{A}>_{E3}$, a single regular integer tessellation of 3-dimensional space into cubes $\{4, 3\}$ is accomplished provided that initial polytope was octahedron $\{3, 4\}$.

If we take a tetrahedron $\{3, 3\}$ as an initial polytope («star of tessellation»), then we will obtain tessellation just with the minimal covering factor = 4. According to (19-20), corresponding regular polytope with $E$-vector [2/3, 1/4] covering 3-dimensional space has the following parameters: $\mathbf{y} = 22.79...$, $\mathbf{x} = 34.19...$, $\mathbf{n} = 13.397...$, $\mathbf{m} = 5.104...$ . It is interesting to note that above values level the average values of Meijering [6], obtained by modeling, with the relative error within 1%. The value $\mathbf{y} = 22.79...$ was first obtained by Smith [12].

Another approach was applied by Donald Coxeter [1]. According to it, regular tessellation of 3-dimensional space (i.e. infinite regular honeycomb) may be regarded as 4-dimensional degenerated regular polytope $\{\mathbf{m}, \mathbf{i}, \mathbf{U}\}$ with infinitely distant center, and meeting the condition:

$$\cos \dfrac{\pi}{\mathbf{i}} = \sin \dfrac{\pi}{\mathbf{m}} \sin \dfrac{\pi}{\mathbf{U}} \; , \quad (21)$$

or the same in the relative bases $[H]$ and $[E]$:

$$\gamma = 1 \; \Rightarrow \; (1-\varepsilon)(1-\eta) - \delta = 0 \; , \quad (22)$$

which has a single integer solution {4, 3, 4} representing an infinite 3-dimensional cubic mosaic. Then, for the reason of closest packing of inscribed spheres, the mosaic {**m**, 3, 3} = [ε, 1/4, 1/4] is considered, where **m**, according to (21), equals to non-integer value **m** = 5,104... (ε = 2/3). This is fully consistent with our results based upon (19-20). In Coxeter's opinion, «non-integer value of **m** means that this «honeycomb» can exist in statistical sense only, though it amazingly agrees with an experimental data».

Transformations $<\mathbf{A}>_{E3}$ make up a cyclic group of 5$^{th}$ order $\mathfrak{C}_5$. That is, consecutive fivefold action of $<\mathbf{A}>_{E3}$ on $E$-vector yields in initial regular polytope irrespectively to which one was such. There are two $E$-vectors invariant to the transformations $<\mathbf{A}>_{E3}$. One of them, with the real values of $f$-symbol {3.47..., 3.47...}, has coordinates $[(3-\sqrt{5})/2, (3-\sqrt{5})/2]$ in a relative basis $[E]$. Hereinafter, $E$-vectors of invariant polytopes will be called eigenvectors of the respective transformations.

It is easy to notice that center of polytope may be reflected not only over vertices, but also over edges and faces, and new polytopes may be obtained similarly. Reflections over edges may be performed in two ways: $\mathbf{i}' = \mathbf{i}$ and $\mathbf{i}' = \mathbf{m}$. The major difference of reflections over edges is that no full tessellation of space by the uniform polytopes is achieved. Instead, each kind of reflection creates its own regular sublattice, and the whole 3-dimensional space appears to be a union of them. Consider both methods of reflections over edges (Fig. 4) which further will be called reflections of the I and II type.

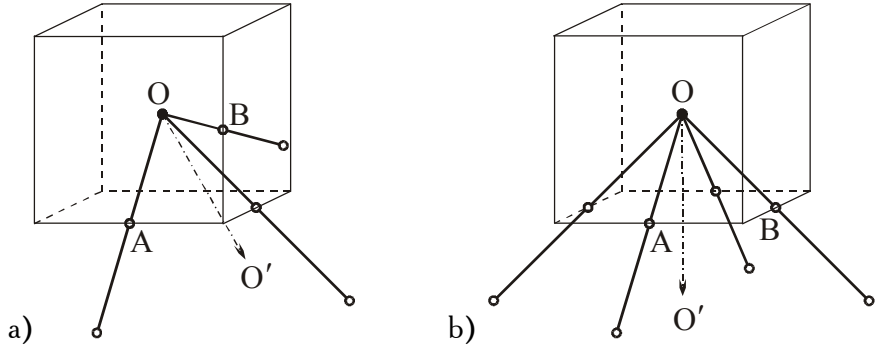

Fig. 4. Reflections of the center of 3-dimensional polytope over edges of a) I and b) II type

Method of calculation of **m**′ is similar to the case of reflection over vertices. Corresponding transformations in the relative basis $[E]$ have the form (we don't show transformations of $f$-symbol here due to its simple relation to the basis $[E]$):

Reflection of I type:

$$<\mathbf{B}>_{E3} = \begin{cases} \varepsilon' = 1 - \varepsilon - \delta \; ; \\ \delta' = \delta \; . \end{cases} \qquad (23)$$

Reflection of II type:

$$<\mathbf{C}>_{E3} = \begin{cases} \varepsilon' = 1 - \varepsilon - \delta \; ; \\ \delta' = \varepsilon \; . \end{cases} \qquad (24)$$

Reflection of I type behaves as involution $\mathfrak{C}_2$, while reflection of II type makes up a cyclic group of the 3$^{rd}$ order $\mathfrak{C}_3$. Together they generate a symmetric group $\mathfrak{S}_3$ of the 6$^{th}$ order. Eigenspace of transformations $<\mathbf{B}>_{E3}$ is constituted by all $E$-vectors satisfying the condition: $\delta = 1 - 2\varepsilon$. In particular, octahedron $\{3, 4\} = [1/4, 1/2]$ has been a single integer regular polytope representing this set. Transformations $<\mathbf{C}>_{E3}$ have a single eigenvector $[1/3, 1/3]$ that is not represented, however, in the set of integer polytopes.

Consider the last method of reflection of the center of 3-dimensional regular polytope - reflection over faces (Fig. 5). In this case, likewise reflecting over vertices, we have a full regular tessellation of space by reflected polytopes, but with a covering factor equal to **y**.

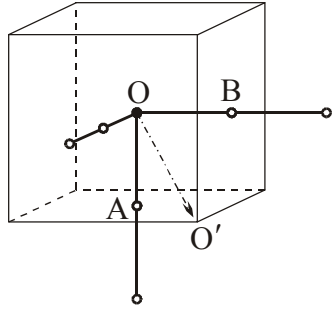

Fig. 5. Reflection of the center of 3-dimensional polytope over faces

Transformations of $E$-vector have a form:

$$<\mathbf{D}>_{E3} = \begin{cases} \varepsilon' = 1 - \dfrac{\delta}{1 - \varepsilon} \; ; \\ \delta' = \delta \; . \end{cases} \qquad (25)$$

Transformations (25) behave as involutions $\mathfrak{C}_2$. All $E$-vectors holding the condition: $\delta = (1 - \varepsilon)^2$ belong to the eigenspace of transformations $<\mathbf{D}>_{E3}$. Among all integer regular polytopes, only $E$-vector of cube satisfies this condition. Thus transformations $<\mathbf{D}>_{E3}$ realize the second (after $<\mathbf{A}>_{E3}$) method of integer regular tessellation of 3-dimensional space by cubes $\{4, 3\}$.

Transformations $<A>_{E3} - <D>_{E3}$ generate an infinite discrete group which will further be called relative 3-dimensional group of reflections of generalized regular polytopes, and be denoted **RRP(3)**. Consider the generating relations in this group (index $E3$ is assumed, compositions are read from left to right, $E$ - identical transformation):

$$<A>^5 = <B>^2 = <C>^3 = <D>^2 = <AC>^4 = <AD>^2 = <BC>^2 = <DCBA> = E, \quad (26)$$

from which, in particular, the following relations may be deduced:

$$<A> = <CBD>; \quad <B> = <ADC>; \quad <C> = <ADB>; \quad <D> = <CBA>. \quad (27)$$

Thus the minimal set of generators of the group **RRP(3)** may be limited to any three elements out from transformations $<A>_{E3} - <D>_{E3}$.

Contemplated group has finite discrete subgroups, in particular already mentioned symmetric subgroup $\mathfrak{S}_3$ of the 6th order generated by the elements $<B>_{E3}$ and $<C>_{E3}$. There are dihedral subgroup $\mathfrak{D}_5$ of the 10th order generated by the elements $<A>_{E3}$ and $<D>_{E3}$, dihedral subgroups $\mathfrak{D}_4$ of the 8th order generated by the sets of elements $\{<AC>_{E3}, <B>_{E3}\}$ or $\{<BC>_{E3}, <D>_{E3}\}$, as well as many various finite and infinite cyclic and other subgroups, determined by the relationships (26) and automorphisms. Finite periods of the **RRP(3)** elements are depleted by the set $\{1, 2, 3, 4, 5\}$.

Note that element $<AD>_{E3}$ provides an involutive transition to the reciprocal polytope:

$$<AD>_{E3} = \begin{cases} \varepsilon' = \delta; \\ \delta' = \varepsilon. \end{cases} \quad (28)$$

From the number of simplest **RRP(3)** transformations, elements $<CDB>_{E3}$, $<DBC>_{E3}$ and $<ACAC>_{E3}$ may be noted as having eigenvectors corresponding to the integer regular polytopes $\{3, 3\}$, $\{4, 3\}$ and $\{3, 4\}$ respectively. It is interesting to note that group **RRP(3)** contains no elements with eigenvector whose $f$-symbol has number 5. In particular, neither dodecahedron $\{5, 3\}$ nor icosahedron $\{3, 5\}$ are the eigenvectors of whichever transformation. The proof of this fact is unknown to the author, though thousands of eigenvectors were tested up to ninefold compositions of transformations $<A>_{E3} - <D>_{E3}$.

Formulation of condition (18) on the existence of regular tessellations presumes consideration of degenerated $(n + 1)$-dimensional regular polytope with infinitely distant center $O_{n+1}$ (or having zero length of an edge). Let $E$-symbol of 5-dimensional regular polytope has components $[\varepsilon, \delta, \eta, \nu]$, and $H$-symbol - components $[\alpha, \beta, \gamma, \mu]$. It is required that

$$\left(^{1}\mathbf{R}_0\right)^2 = \left(^{5}\mathbf{R}_0\right)^2 - \left(^{5}\mathbf{R}_1\right)^2 = 0 \implies \mu \equiv \frac{\left(^{5}\mathbf{R}_0\right)^2}{\left(^{5}\mathbf{R}_1\right)^2} = 1, \tag{29}$$

where $\mu$, in concordance to (6-9), is related to the components of $E$-vector as follows:

$$\mu = \frac{1}{\varepsilon}\left(1 - \frac{\delta(1-\nu)}{1-\eta-\nu}\right), \tag{30}$$

which, if $\mu = 1$, yields to:

$$\nu = 1 - \frac{\eta(1-\varepsilon)}{1-\varepsilon-\delta}. \tag{31}$$

Thus $E$-vector of polytope $[\delta, \eta, \nu]$ (i.e. «star of tessellation») may be obtained from the initial polytope $[\varepsilon, \delta, \eta]$ by the transformation:

$$\begin{aligned}\varepsilon' &= \delta \,; \\ \delta' &= \eta \,; \\ \eta' &= 1 - \frac{\eta(1-\varepsilon)}{1-\varepsilon-\delta} \,.\end{aligned} \tag{32}$$

It is easy to check that all known regular tessellations of 4-dimensional space by integer regular polytopes satisfy the above relationships. Similarly, transformations for any $n$ may be obtained.

Adjusting formulation of Berger [2] (18), we note that existence of polytope $\{f_2,..., f_n\} = [\varepsilon', \delta', \eta']$ is limited just by the condition $\varepsilon + \delta \neq 1$, and only its «integerness» makes a difference. Yet, according to (32), in order to check existence of regular tessellation, there is no need to consider degenerated $(n+1)$-dimensional polytopes.

Further we will show that condition (32) of the existence of regular tessellations of 4-dimensional space is equivalent to the condition of existence of regular polytope obtained as a result of inverse transformation of $E$-symbol when reflecting polytope over vertices, and that such transformation does not deplete all possible methods of regular tessellation.

### Systematic construction of tessellations of 4-dimensional space

Tessellation of space into regular honeycombs is performed by reflections of the center of regular polytope over the centers of its $(0, ..., n-1)$-dimensional elements («faces»). It is sufficient to construct $n$ linearly independent constituting vectors of reflected polytope.

Unfortunately, reflections of 4-dimensional polytope have not been as illustrative as those in the case of 3 dimensions.

Let introduce complementary 3-index vectors ${}^{ik}\mathbf{r}_j$ in hyperplanes passing through centers of $i$, $j$ and $k$-faces ($0 \leq i < j < k \leq n$):

$$ {}^{ik}\mathbf{r}_j = {}^{k}\mathbf{R}_j - \frac{\left({}^{k}\mathbf{R}_j\right)^2}{\left({}^{k}\mathbf{R}_i\right)^2} {}^{k}\mathbf{R}_i \ ; \tag{33}$$

$$ \left({}^{ik}\mathbf{r}_j\right)^2 = \frac{\left({}^{j}\mathbf{R}_i\right)^2 \left({}^{k}\mathbf{R}_j\right)^2}{\left({}^{k}\mathbf{R}_i\right)^2} \ . \tag{34}$$

**1. Reflection over 0-faces (vertices)**

$$<\mathbf{A}>_4 \ = \ \begin{cases} {}^{1}\mathbf{R}'_0 = -\mathbf{p}_0 \\ {}^{02}\mathbf{r}'_1 = {}^{1}\mathbf{R}_0 \\ {}^{03}\mathbf{r}'_1 = {}^{2}\mathbf{R}_0 \\ {}^{04}\mathbf{r}'_1 = {}^{3}\mathbf{R}_0 \end{cases} \tag{35}$$

$$<\mathbf{A}>_{E4} \ = \ \begin{cases} \varepsilon' = 1 - \dfrac{\varepsilon(1-\eta)}{1-\delta-\eta} \\ \delta' = \varepsilon \\ \eta' = \delta \end{cases} \tag{36}$$

$$<\mathbf{A}>_{P4} \ = \ \begin{cases} \mathbf{p}'_0 = -\dfrac{(1-\varepsilon-\delta)(1-\delta-\eta)}{\varepsilon\,\delta\,\eta} \mathbf{p}_3 \\ \mathbf{p}'_1 = \mathbf{p}_0 + \mathbf{p}'_0 \\ \mathbf{p}'_2 = \dfrac{1-\delta-\eta}{\varepsilon(1-\eta)} \mathbf{p}_1 + \mathbf{p}'_0 \\ \mathbf{p}'_3 = \dfrac{(1-\varepsilon)(1-\delta-\eta)}{\varepsilon\,\delta} \mathbf{p}_2 + \mathbf{p}'_0 \end{cases} \tag{37}$$

Transformations $<\mathbf{A}>_{E4}$ make up a cyclic group of the 6$^{\text{th}}$ order $\mathfrak{C}_6$. Eigenvectors of these transformations are $E$-vectors [1/3, 1/3, 1/3] and [1, 1, 1]. Reflections over 0-faces (vertices) provide a full regular tessellation of 4-dimensional space by regular polytopes.

Yet, it should be noted that transformations (32), derived earlier as a condition of existence of regular tessellation of 4-dimensional space, are in fact inverse transformations with regard to reflections over vertices $<\mathbf{A}>_{E4}$, that can easily be checked by direct calculations.

**2. Reflection over 1-faces (edges) of I type**

$$<\mathbf{B}>_4 \;=\; \begin{cases} {}^1\mathbf{R}'_0 = -\mathbf{p}_1 \\ {}^{02}\mathbf{r}'_1 = {}^{02}\mathbf{r}_1 \\ {}^{03}\mathbf{r}'_1 = {}^{03}\mathbf{r}_1 \\ {}^{04}\mathbf{r}'_1 = {}^{04}\mathbf{r}_1 \end{cases} \tag{38}$$

$$<\mathbf{B}>_{E4} \;=\; \begin{cases} \varepsilon' = 1 - \varepsilon - \dfrac{\delta}{1-\eta} \\ \delta' = \delta \\ \eta' = \eta \end{cases} \tag{39}$$

$$<\mathbf{B}>_{P4} \;=\; \begin{cases} \mathbf{p}'_0 = -\mathbf{p}_0 \\ \mathbf{p}'_1 = -\mathbf{p}_0 + \mathbf{p}_1 \\ \mathbf{p}'_2 = \dfrac{1}{\delta + \varepsilon(1-\eta)} \left( -\delta\,\mathbf{p}_0 + (1-\varepsilon)(1-\eta)\,\mathbf{p}_2 \right) \\ \mathbf{p}'_3 = \dfrac{1}{\varepsilon + \eta(\delta - \varepsilon)} \left( -\delta\eta\,\mathbf{p}_0 + (1-\varepsilon-\eta)(1-\eta)\,\mathbf{p}_3 \right) \end{cases} \tag{40}$$

Transformations $<\mathbf{B}>_{E4}$ behave as involutions $\mathfrak{C}_2$. Eigenspace of $<\mathbf{B}>_{E4}$ is constituted by the set of all $E$-vectors satisfying the condition: $\delta = (1 - 2\varepsilon)(1 - \eta)$.

### 3. Reflection over 1-faces (edges) of II type

$$<\mathbf{C}>_4 \;=\; \begin{cases} {}^1\mathbf{R}'_0 = -\mathbf{p}_1 \\ {}^{02}\mathbf{r}'_1 = {}^{02}\mathbf{r}_1 \\ {}^{03}\mathbf{r}'_1 = {}^{03}\mathbf{r}_1 \\ {}^{04}\mathbf{r}'_1 = {}^3\mathbf{R}_1 \end{cases} \tag{41}$$

$$<\mathbf{C}>_{E4} \;=\; \begin{cases} \varepsilon' = 1 - \varepsilon - \dfrac{\delta}{1-\eta} \\ \delta' = \delta \\ \eta' = \dfrac{\varepsilon}{\varepsilon + \delta} \end{cases} \tag{42}$$

$$<\mathbf{C}>_{P4} \;=\; \begin{cases} \mathbf{p}'_0 = -\dfrac{(1-\varepsilon-\delta)(1-\eta)}{\delta\,\eta}\,\mathbf{p}_3 \\ \mathbf{p}'_1 = \mathbf{p}_1 + \mathbf{p}'_0 \\ \mathbf{p}'_2 = \dfrac{1}{\delta + \varepsilon(1-\eta)} \left( \varepsilon(1-\eta)\,\mathbf{p}_0 + (1-\varepsilon)(1-\eta)\,\mathbf{p}_2 \right) + \mathbf{p}'_0 \\ \mathbf{p}'_3 = \dfrac{\varepsilon(1-\eta)}{\varepsilon + \eta(\delta-\varepsilon)} \left( \mathbf{p}_0 + \mathbf{p}'_0 \right) \end{cases} \tag{43}$$

Transformations $<C>_{E4}$ make up a cyclic group of the 3$^{rd}$ order $\mathfrak{C}_3$. Their eigenspace is constituted by the set of all *E*-vectors satisfying the condition: $\{\ \varepsilon = \eta/(1+2\eta)\ ,\ \delta = (1-\eta)/(1+2\eta)\ \}$.

### 4. Reflection over 1-faces (edges) of III type

$$<D>_4 = \begin{cases} {}^1R'_0 = -p_1 \\ {}^{02}r'_1 = {}^{02}r_1 \\ {}^{03}r'_1 = {}^2R_1 \\ {}^{04}r'_1 = {}^3R_1 \end{cases} \tag{44}$$

$$<D>_{E4} = \begin{cases} \varepsilon' = 1 - \varepsilon - \dfrac{\delta}{1-\eta} \\ \delta' = \varepsilon \\ \eta' = \dfrac{\delta}{\varepsilon + \delta} \end{cases} \tag{45}$$

$$<D>_{P4} = \begin{cases} p'_0 = -\dfrac{(1-\varepsilon-\delta)(1-\eta)}{\delta\,\eta}\,p_3 \\ p'_1 = p_1 + p'_0 \\ p'_2 = \dfrac{1}{\delta + \varepsilon(1-\eta)}\bigl(\varepsilon(1-\eta)\,p_0 + (1-\varepsilon)(1-\eta)\,p_2\bigr) + p'_0 \\ p'_3 = \dfrac{(1-\varepsilon)(1-\eta)}{\delta}\,p_2 + p'_0 \end{cases} \tag{46}$$

Transformations $<D>_{E4}$ make up a cyclic group of the 4$^{th}$ order $\mathfrak{C}_4$. They have a single eigenvector [1/4, 1/4, 1/2] corresponding to the integer regular polytope {3, 3, 4}. That is, this polytope may be «self-reproduced» with the transformations $<D>_{E4}$.

Reflections of polytopes over 1-faces of I, II, and III types generate together a symmetric group $\mathfrak{S}_4$ of the 24$^{th}$ order. Likewise the case of 3 dimensions, reflections over edges does not provide a full regular tessellation of 4-dimensional space individually. However, the union of 3 regular sublattices, corresponding to the reflections of I, II, and III types, fully covers 4-dimensional space.

### 5. Reflection over 2-faces (faces) of I type

$$<E>_4 = \begin{cases} {}^1R'_0 = -p_2 \\ {}^{02}r'_1 = {}^{13}r_2 \\ {}^{03}r'_1 = {}^{14}r_2 \\ {}^{04}r'_1 = {}^{04}r_2 \end{cases} \tag{47}$$

$$<\mathbf{E}>_{E4} = \begin{cases} \varepsilon' = 1 - \eta - \dfrac{\delta}{1-\varepsilon} \\ \delta' = \eta \\ \eta' = \dfrac{\delta}{\delta + \eta} \end{cases} \tag{48}$$

$$<\mathbf{E}>_{P4} = \begin{cases} \mathbf{p}'_0 = -\mathbf{p}_0 \\ \mathbf{p}'_1 = -\mathbf{p}_0 + \mathbf{p}_2 \\ \mathbf{p}'_2 = -\mathbf{p}_0 + \dfrac{1}{\delta + \eta(1-\varepsilon)}\left(\delta\,\mathbf{p}_1 + (1-\varepsilon-\delta)\,\mathbf{p}_3\right) \\ \mathbf{p}'_3 = -\mathbf{p}_0 + \mathbf{p}_1 \end{cases} \tag{49}$$

Unlike all other reflections mentioned so far, transformations $<\mathbf{E}>_{E4}$ don't have a finite period. *E*-vectors [0, 1/2, 1/2] and [3/2, 1/2, 1/2] are the eigenvectors of these transformations.

### 6. Reflection over 2-faces (faces) of II type

$$<\mathbf{F}>_4 = \begin{cases} {}^1\mathbf{R}'_0 = -\mathbf{p}_2 \\ {}^{02}\mathbf{r}'_1 = {}^{13}\mathbf{r}_2 \\ {}^{03}\mathbf{r}'_1 = {}^{03}\mathbf{r}_2 \\ {}^{04}\mathbf{r}'_1 = {}^{04}r_2 \end{cases} \tag{50}$$

$$<\mathbf{F}>_{E4} = \begin{cases} \varepsilon' = 1 - \eta - \dfrac{\delta}{1-\varepsilon} \\ \delta' = \delta \\ \eta' = \dfrac{\eta}{\delta + \eta} \end{cases} \tag{51}$$

$$<\mathbf{F}>_{P4} = \begin{cases} \mathbf{p}'_0 = -\mathbf{p}_0 \\ \mathbf{p}'_1 = -\mathbf{p}_0 + \mathbf{p}_2 \\ \mathbf{p}'_2 = -\mathbf{p}_0 + \dfrac{1}{\delta + \eta(1-\varepsilon)}\left(\delta\,\mathbf{p}_1 + (1-\varepsilon-\delta)\,\mathbf{p}_3\right) \\ \mathbf{p}'_3 = \dfrac{1}{\eta + \varepsilon(\delta - \eta)}\left(-\eta(1-\varepsilon)\,\mathbf{p}_0 + (1-\varepsilon-\delta)\,\mathbf{p}_3\right) \end{cases} \tag{52}$$

Transformations $<\mathbf{F}>_{E4}$ don't have a finite period. Eigenspace of $<\mathbf{F}>_{E4}$ is constituted by three sets of all *E*-vectors satisfying the conditions: $\{\eta = 0,\ \delta = (1-\varepsilon)^2\}$, $\{\varepsilon = 0,\ \delta = 1 - \eta\}$ and $\{\varepsilon = 2 - \eta,\ \delta = 1 - \eta\}$.

### 6. Reflection over 2-faces (faces) of III type

$$<\mathbf{G}>_4 \;=\; \begin{cases} {}^1\mathbf{R}'_0 = -\mathbf{p}_2 \\ {}^{02}\mathbf{r}'_1 = {}^{13}\mathbf{r}_2 \\ {}^{03}\mathbf{r}'_1 = {}^{03}\mathbf{r}_2 \\ {}^{04}\mathbf{r}'_1 = {}^{3}\mathbf{R}_2 \end{cases} \tag{53}$$

$$<\mathbf{G}>_{E4} \;=\; \begin{cases} \varepsilon' = 1 - \eta - \dfrac{\delta}{1-\varepsilon} \\ \delta' = \delta \\ \eta' = \varepsilon \end{cases} \tag{54}$$

$$<\mathbf{G}>_{P4} \;=\; \begin{cases} \mathbf{p}'_0 = -\dfrac{1-\varepsilon-\delta}{\eta(1-\varepsilon)}\mathbf{p}_3 \\ \mathbf{p}'_1 = \mathbf{p}_2 + \mathbf{p}'_0 \\ \mathbf{p}'_2 = \dfrac{\delta}{\delta + \eta(1-\varepsilon)}(\mathbf{p}_1 + \mathbf{p}'_0) \\ \mathbf{p}'_3 = \dfrac{\varepsilon\delta}{\eta + \varepsilon(\delta-\eta)}(\mathbf{p}_0 + \mathbf{p}'_0) \end{cases} \tag{55}$$

Likewise $<\mathbf{E}>_{E4}$ and $<\mathbf{F}>_{E4}$, transformations $<\mathbf{G}>_{E4}$ don't have a finite period. Eigenspace of $<\mathbf{G}>_{E4}$ is constituted by the set of all *E*-vectors satisfying the condition: $\{\eta = \varepsilon,\; \delta = (1-\varepsilon)(1-2\varepsilon)\}$.

Reflections over 2-faces of I, II and III types don't provide a full regular tessellation of 4-dimensional space individually. Each of them creates its own regular sublattice, and their union covers the whole 4-dimensional space.

### 7. Reflection over 3-faces (3-dimensional polytopes)

$$<\mathbf{H}>_4 \;=\; \begin{cases} {}^1\mathbf{R}'_0 = -\mathbf{p}_3 \\ {}^{02}\mathbf{r}'_1 = {}^{24}\mathbf{r}_3 \\ {}^{03}\mathbf{r}'_1 = {}^{14}\mathbf{r}_3 \\ {}^{04}\mathbf{r}'_1 = {}^{04}\mathbf{r}_3 \end{cases} \tag{56}$$

$$<\mathbf{H}>_{E4} \;=\; \begin{cases} \varepsilon' = 1 - \dfrac{\eta(1-\varepsilon)}{1-\varepsilon-\delta} \\ \delta' = \eta \\ \eta' = \delta \end{cases} \tag{57}$$

$$<\mathbf{H}>_{P4} \;=\; \begin{cases} \mathbf{p}'_0 = -\mathbf{p}_0 \\ \mathbf{p}'_1 = -\mathbf{p}_0 + \mathbf{p}_3 \\ \mathbf{p}'_2 = -\mathbf{p}_0 + \mathbf{p}_2 \\ \mathbf{p}'_3 = -\mathbf{p}_0 + \mathbf{p}_1 \end{cases} \tag{58}$$

Transformations $<H>_{E4}$ behave as involutions $\mathfrak{C}_2$. Eigenspace of $<H>_{E4}$ is constituted by two sets of all *E*-vectors satisfying the conditions: { $\varepsilon = 1$, $\delta = \eta$ } and { $\varepsilon = 1 - 2\eta$, $\delta = \eta$ }. It should be noted that transformations $<H>_{E4}$ realize the second (after $<A>_{E4}$) method of a full regular tessellation of 4-dimensional space.

Transformations $<A>_{E4} - <H>_{E4}$ generate an infinite discrete group which will further be called relative 4-dimensional group of reflections of generalized regular polytopes, and be denoted **RRP(4)**. Consider the generating relations in this group (index *E4* is assumed, compositions are read from left to right, small letters denote inverse elements, *E* - identical transformation):

$$<A>^6 = <B>^2 = <C>^3 = <D>^4 = <H>^2 = <AE>^2 = <AF>^3 = <AG>^4 = <AH>^2 =$$

$$<BC>^2 = <BD>^3 = <BG>^2 = <CD>^4 = <cF>^2 = <De>^2 = <Ef>^2 = <Eg>^3 = <Fg>^2 =$$

$$<AAD>^{10} = <AAG>^4 = <ACe>^6 = <ACg>^{10} = <AEE>^6 = <AEH>^{10} = <BHC>^{10} =$$

$$<AECB> = <AHDe> = <BAHG> = <AHCf> = <HGDDg> = \ldots = E, \quad (59)$$

from which, in particular, the following relations may be deduced:

$$<B> = <AHG>; \quad <C> = <dBd^2>; \quad <E> = <AHD>; \quad <F> = <AHC>; \quad <H> = <GDDg>.$$
$$(60)$$

Thus the minimal set of generators of the group **RRP(4)** may be limited to three elements. For example, it may be represented by transformations $<A>_{E4}$, $<D>_{E4}$, and $<G>_{E4}$.

As it can be seen, the contemplated group is much more complicated than **RRP(3)**. Reflection over vertices $<A>_{E4}$ and over 3-dimensional polytopes $<H>_{E4}$ generate together dihedral subgroup $\mathfrak{D}_6$ of $12^{th}$ order. Element of this subgroup $<AH>_{E4}$ having a period 2 is remarkable as it provides an involutive transition to the reciprocal polytope:

$$<AH>_{E4} = \begin{cases} \varepsilon' = \eta; \\ \delta' = \delta; \\ \eta' = \varepsilon. \end{cases} \quad (61)$$

Group **RRP(4)** contains dihedral subgroup $\mathfrak{D}_{10}$ of the $20^{th}$ order generated by the elements $<AE>_{E4} = <CB>_{E4}$ and $<H>_{E4}$. The relation $<AEH>^{10}_{E4} = E$ may be checked directly, though the manual calculations are rather cumbersome. Mentioned earlier the symmetric subgroup $\mathfrak{S}_4$ of the $24^{th}$ order is generated by the sets of elements {$<B>_{E4}$, $<D>_{E4}$} or {$<C>_{E4}$, $<D>_{E4}$}. Group **RRP(4)** has also many various finite and infinite cyclic and other subgroups determined by relationships (59) and automorphisms. Finite periods of **RRP(4)** elements are depleted by the set {1, 2, 3, 4, 5, 6, 10}.

Out from the number of simplest transformations of **RRP(4)**, elements $<\mathbf{D}>_{E4}$, $<\mathbf{EAH}>_{E4}$ and $<\mathbf{HAH}>_{E4}$ may be noted as having eigenvectors equal to *E*-vectors of integer regular polytopes {3, 3, 4}, {4, 3, 3} and {3, 4, 3} respectively. Likewise **RRP(3)**, the group **RRP(4)** apparently contains no elements with an eigenvector corresponding to *f*-symbol with the digit 5. In particular, neither Euclidean polytopes {5, 3, 3} and {3, 3, 5}, nor pseudo-Euclidean polytopes {5, 3, 4}, {4, 3, 5}, {3, 5, 3} and {5, 3, 5} are the eigenvectors of whichever transformations. Yet, we couldn't find any element in **RRP(4)** having a *single* eigenvector corresponding to the polytope {3, 3, 3}.

We shall assign special term **eigentopes** to the self-reproducing polytopes with finite eigenvectors.

### Isotropy subgroups of the group RRP(4) and their linear representation

**Definition.** Isotropy subgroup of an arbitrary *E*-vector $x$ (and corresponding regular polytope) is a set $\mathbf{Iso}(x) \subset \mathbf{RRP(4)}$: $\{g \in \mathbf{RRP(4)} : g(x) = x\}$.

As it follows from the above definition, all polytopes whose *E*-vectors belong to the eigenspaces (and in particular, eigenvectors) of transformations **RRP(4)**, have non-trivial isotropy subgroups. For example, (infinite) isotropy subgroup $\mathbf{Iso}(\{3, 3, 4\})$ of the polytope {3, 3, 4} is generated by the elements $<\mathbf{A}^3>_{E4}$, $<\mathbf{B}>_{E4}$ and $<\mathbf{D}>_{E4}$. Another integer Euclidean and pseudo-Euclidean 4-dimensional polytopes also have non-trivial isotropy subgroups.

For each isotropy subgroup $\mathbf{Iso}(x)$ of the group **RRP(4)**, a linear representation into the matrix group $\mathbf{Q}(x) \subset GL(3, \mathbf{R})$ may be constructed. Such representation appears to be strict for the finite isotropy subgroups.

### Absolute group of reflections of 4-dimensional regular polytopes

Group **RRP(4)** acts in the space of *E*-vectors of generalized regular polytopes (**RRP(4)**-module) irrespectively to change of linear dimensions of reflected polytopes and orientation of their constituting frames. As a matter of fact, identical transformations of an arbitrary *E*-vector of the form $<\mathbf{RST}...>_{E4}^{q} = E$ as well as transformations of *E*-eigenvector under the action of its isotropy subgroup, in general result in change of scale (and therefore scalar squares of the linear elements) of polytopes by a certain factor. Yet, in general, there is also change in orientation of constituting frame of a polytope.

As an example, consider reflections over vertices (36-37). Unlike the transformations of *E*-vectors in a relative basis for which the relationship $<\mathbf{A}>^6_{E4} = E$ is satisfied, sixfold transformations of an arbitrary **P**-vector with the frame {**P**} under the action of $<\mathbf{A}>_{P4}$ have not been identical:

$$<\mathbf{A}>^6_{P4} = \begin{cases} \mathbf{p}_0^{VI} = \lambda \mathbf{p}_0 \\ \mathbf{p}_1^{VI} = \lambda \mathbf{p}_1 \\ \mathbf{p}_2^{VI} = \lambda \mathbf{p}_2 \\ \mathbf{p}_3^{VI} = \lambda \mathbf{p}_3 \end{cases} ; \quad \lambda = -\frac{(1-\varepsilon-\delta)(1-\delta-\eta)}{\varepsilon\delta\eta(1-\varepsilon-\delta-\eta+\varepsilon\eta)}. \tag{62}$$

One of the eigenvectors of transformation $<\mathbf{A}>_{E4}$ [1/3, 1/3, 1/3] corresponds to the eigenvalue $\lambda = -27$ of the transformation $<\mathbf{A}>^6_{P4}$. Thus sixfold action of $<\mathbf{A}>_{P4}$ on the eigentope [1/3, 1/3, 1/3] causes its scale to increase by 27 times, and the frame {**P**} orientation to inverse. A single action of $<\mathbf{A}>_{P4}$ causes scale of this polytope to increase by $|\lambda|^{1/6} = \sqrt{3}$ times; and it is accompanied by the discrete rotation of transformed frame so that $\cos(\mathbf{P}', \mathbf{P}) = -1/\sqrt{3}$.

Unlike [1/3, 1/3, 1/3], sixfold action of transformation $<\mathbf{A}>_{P4}$ on another eigenvector of $<\mathbf{A}>_{E4}$ [1, 1, 1], associated with an eigenvalue $\lambda = 1$, causes no changes neither in the scale of the (degenerated) eigentope, nor in the orientation of its frame {**P**}.

Thus absolute transformations $<\mathbf{A}>_{P4} - <\mathbf{H}>_{P4}$ (37-58) carry an additional information on change in linear dimensions of transformed polytopes, and the orientation of their constituting frames. These transformations generate an infinite discrete group, which we will call an absolute group of reflections of 4-dimensional generalized regular polytopes, and denote it as **ARP(4).** Consider the generating relations in this group (index **P4** is assumed, compositions are read from left to right, small letters denote inverse elements, *E* - identical transformation):

$$<\mathbf{B}>^2 = <\mathbf{C}>^3 = <\mathbf{D}>^4 = <\mathbf{H}>^2 = <\mathbf{BC}>^2 = <\mathbf{BD}>^3 = <\mathbf{CD}>^4 = <\mathbf{Ef}>^2 = <\mathbf{Eg}>^3 = <\mathbf{Fg}>^2 =$$
$$<\mathbf{BHC}>^{10} = <\mathbf{BDeG}> = <\mathbf{HECg}> = <\mathbf{AE}>^2<\mathbf{AH}>^{-2} = <\mathbf{AGEHg}>^4 = <\mathbf{GBCf}> =$$
$$<\mathbf{egDaGH}> = <\mathbf{HGDDg}> = \ldots = E, \tag{63}$$

from which, in particular, the following relations may be deduced:

$$<\mathbf{B}> = <\mathbf{DCD}^2> ; \ <\mathbf{C}> = <\mathbf{eHG}> ; \ <\mathbf{E}> = <\mathbf{gDaGH}>; \ <\mathbf{F}> = <\mathbf{GBC}>; \ <\mathbf{H}> = <\mathbf{GDDg}>. \tag{64}$$

Thus, likewise **RRP(4)**, the minimal set of generators of the group **ARP(4)** may be limited to three elements, and represented e. g. by the following transformations: $<\mathbf{A}>_{P4}$, $<\mathbf{D}>_{P4}$, and $<\mathbf{G}>_{P4}$.

Group **RRP(4)** has been a factor-group of the group **ARP(4)**. An infinite invariant subgroup **NRP(4)** $\subset$ **ARP(4)** appears to be a kernel of homomorphism with the elements:

$$<\mathbf{A}>^{6k},\ <\mathbf{AE}>^{2k},\ <\mathbf{AG}>^{4k},\ <\mathbf{AH}>^{2k},\ <\mathbf{BG}>^{2k},\ <\mathbf{AAG}>^{4k},\ <\mathbf{ACe}>^{6k},\ <\mathbf{ACg}>^{10k},$$

$$<\mathbf{AEE}>^{6k},\ <\mathbf{AEH}>^{10k},\ <\mathbf{AECB}>^{k},\ \ldots,\ \ k \in \mathbf{Z}$$

being mapped to the unit element of the group **RRP(4)**.

### Matrix representation of the group ARP(4)

Absolute transformations (37-58) may be expressed in form of matrix operators of polytope reflections **W**, multiplicatively acting on the natural frame {**P**}:

$$\mathbf{P}' = \mathbf{W}\,\mathbf{P}, \tag{65}$$

where rows of matrices **P** and **P**′ contain components of vectors $\{\mathbf{p}_0, \mathbf{p}_1, \mathbf{p}_2, \mathbf{p}_3\}$ of constituting frame of initial and transformed polytopes respectively. Basis vectors are the vectors of the frame of the initial polytope. Hence **P** is an identity matrix 4×4.

The following matrix, for example, represents the transformation $<\mathbf{B}>_{\mathbf{P4}}$ :

$$\mathbf{W}(<\mathbf{B}>_{\mathbf{P4}}) = \begin{pmatrix} -1 & 0 & 0 & 0 \\ -1 & 1 & 0 & 0 \\ -\dfrac{\delta}{\delta + \varepsilon(1-\eta)} & 0 & \dfrac{(1-\varepsilon)(1-\eta)}{\delta + \varepsilon(1-\eta)} & 0 \\ -\dfrac{\delta\eta}{\varepsilon + \eta(\delta-\varepsilon)} & 0 & 0 & \dfrac{(1-\varepsilon-\eta)(1-\eta)}{\varepsilon + \eta(\delta-\varepsilon)} \end{pmatrix} \tag{66}$$

In general case, elements of matrix operators **W** are functions of *E*-vector components of polytope, therefore their consecutive action means in general inequality of matrices **W**′ and **W** corresponding to the same transformations:

$$\mathbf{P}'' = \mathbf{W}'\mathbf{P}' = \mathbf{W}'\mathbf{W}\,\mathbf{P}, \tag{67}$$

since *E*-vector, provided that it does not belong to the eigenspace of **W**, is transformed under the action of operators $<>_{E4}$.

For the eigentopes, metric tensors are transformed as follows:

$$G' = \lambda G = \mathbf{W} G \mathbf{W}^T, \tag{68}$$

where

$$\lambda^2 = |\det \mathbf{W}|. \tag{69}$$

This means that transformations of the frames of eigentopes are in fact being *conformal*.

### Spin properties of eigentopes (self-reproducing polytopes)

Consider $E$-eigenvector for the RRP(4) transformation $\langle \mathbf{RST...}\rangle_{E4}$. How do the corresponding ARP(4) conformal transformations $\mathbf{X} = \langle \mathbf{RST...}\rangle_{P4}$ act on the frame of the eigentope? Is there finite value $q$ such that $\mathbf{X}^q = \lambda \mathbf{Id}$ or $\mathbf{U}^q = \mathbf{Id}$ where $\mathbf{U} = (1/\lambda)^{1/q}\mathbf{X}$ is a (pseudo-) orthogonal matrix. We already know the answer for the eigenvectors of ARP(4) as stated above (62).

Computer modelling with a great variety of up to 9-fold $E$-eigenvectors in pseudo-Euclidean space shows that sought values of $q$ do exist:

$q = 2, 3, 4, 5, 6, 7, 8, 9, 10, 12, 13, 14, 16$ for $\lambda > 0$,

and $q = 3, 5, 7, 9, 13$ for $\lambda < 0$.

Transformations $\mathbf{U}$ constitute pseudo-orthogonal cyclic discrete groups $G_U \subset O(1,3)$ with respect to their eigenvectors. The order of these groups $q$ can be considered as an internal degree of freedom of quantum microobjects that can be associated with eigentopes. We shall define the **spin** property of such eigentopes as $J = (q-1)/2$.

The simplest pseudo-Euclidean eigentopes with the finite *spin* values are shown in the following table:

| $q$ | $J$ | Word | $\varepsilon$ | $\delta$ | $\eta$ | $\lambda > 0$ | Word | $\varepsilon$ | $\delta$ | $\eta$ | $\lambda < 0$ |
|---|---|---|---|---|---|---|---|---|---|---|---|
| 2 | 1/2 | agdfcF | 0.577 | 0.366 | 0.366 | 0.179 | | | | | |
| 3 | 1 | ABaGE | 0.382 | 0.382 | 0.500 | 0.056 | ABE$^2$F | 0.293 | 0.414 | 0.547 | -0.071 |
| 4 | 3/2 | AGA | 0.500 | 0.281 | 0.500 | 7.507 | | | | | |
| 5 | 2 | ADagdE | 0.236 | 0.382 | 0.567 | 1.000 | A$^2$FafeC | 0.236 | 0.500 | 0.433 | -0.455 |
| 6 | 5/2 | ADGFAC | 0.360 | 0.360 | 0.500 | 0.227 | | | | | |
| 7 | 3 | aE$^2$HD | 0.298 | 0.537 | 0.409 | 0.213 | aE$^2$aE | 0.298 | 0.537 | 0.409 | -61.7 |
| 8 | 7/2 | HAeda | 0.382 | 0.553 | 0.276 | 0.341 | | | | | |
| 9 | 4 | ABGDFHCF | 0.420 | 0.210 | 0.756 | 3.1E-9 | ABGDFaFF | 0.420 | 0.210 | 0.756 | -1.000 |
| 10 | 9/2 | aEGAdg | 0.691 | 0.236 | 0.618 | 0.120 | | | | | |
| 12 | 11/2 | CAeCAgd | 0.366 | 0.500 | 0.446 | 4096 | | | | | |
| 13 | 6 | aFEGaFEF | 0.652 | 0.250 | 0.705 | 1696 | aFEGHCEF | 0.652 | 0.250 | 0.705 | -5.099 |
| 14 | 13/2 | aHgefgd | 0.555 | 0.308 | 0.555 | 173.1 | | | | | |
| 16 | 15/2 | aeDHDgd | 0.707 | 0.185 | 0.631 | 72.1 | | | | | |

**Discussion on the obtained results**

Perhaps, the most exciting applications of generalized regular polytopes can be suggested in physics of elementary particles, whose behavior and transformations resemble many features of eigentopes. The following analogies can be considered:

| | |
|---|---|
| Infinite number of eigentopes with metric $(+---)$. | Infinite number of elementary particles. |
| Discrete transformations in $ARP(4)$. | Discrete nature of mutual transformations of elementary particles. |
| Statistical interpretation of eigentopes. | Fundamental uncertainty in a microworld. |
| Pseudo-orthogonal dicrete group of $q^{th}$ order of conformal frame transformations of an eigentope. | Spin of an elementary particle $J = (q-1)/2$. |
| Expansion (shrinking) of a frame of an eigentope. | Mass of an elementary particle. |

In the author's opinion, the following reflections must be also taken into consideration:

1. Eigentopes do not require the space-time to have a certain metric. Instead, they determine the metric themselves, likewise energy-momentum tensor does.

2. Eigentope's frame orientation has no preferential direction. This means that observation of a single eigentope object does not depend on a viewing direction. As a result, projection of a spin property on any coordinate axis will be identical. In addition, multiple eigentope objects with the same $E$-eigenvector will be indistinguishable.

3. Any act of observation of eigentope object requires synchronization of its frame with the reference frame of the observer (e.g. the probe). Otherwise eigentope object is invisible. Such observation breaks statistical state of an eigentope and forces one of possible discrete states to transform to.

4. Despite that we know everything about statistical state of an eigentope (*real* Schlafli symbol, frame vectors, and transformation matrix of the $q$-th order), there is still no way to relate this state with a physical field.

**Conclusions**

Generalized regular polytopes can be considered as a new class of geometrical objects. All areas of knowledge where regular polytopes are applied now can absorb new ideas about generalized regular polytopes.

This work does not pretend on a complete theory of generalized regular polytopes and their transformations. It is just an invitation to open wide discussion and initiate further studies on this promising topic.

26